\documentclass{amsart}
\usepackage{amssymb,amscd,amsmath,eucal}
\usepackage[matrix,arrow,curve,frame]{xy}    % XY-pic diagram pac
\xymatrixcolsep{1.9pc}                          % Adjust size of diagrams.
\xymatrixrowsep{1.9pc}
\newdir{ >}{{}*!/-8pt/\dir{>}}                  % Make better tailed arrows

\makeatletter
\renewcommand{\phi}{\varphi}
\renewcommand{\geq}{\geqslant}
\renewcommand{\leq}{\leqslant}
\renewcommand{\epsilon}{\varepsilon}

\renewcommand{\kappa}{\varkappa}

\DeclareMathOperator{\spec}{\sf Spec} \DeclareMathOperator{\ann}{\sf
ann} 
\DeclareMathOperator{\serre}{Serre} \DeclareMathOperator{\inj}{\sf
Inj}  \DeclareMathOperator{\supp}{\sf
supp} \DeclareMathOperator{\open}{{open}}
\DeclareMathOperator{\injzar}{\sf Inj_{\rm zar}}
\DeclareMathOperator{\injzg}{\sf Inj_{\rm zg}}
\DeclareMathOperator{\qinjzar}{\sf QInj_{\rm zar}}
\DeclareMathOperator{\qinjzg}{\sf QInj_{\rm zg}}
\DeclareMathOperator{\End}{End} 
 
 \DeclareMathOperator{\Inj}{\sf Inj}
\DeclareMathOperator{\QInj}{\sf QInj} \DeclareMathOperator{\Gr}{Gr}
 \DeclareMathOperator{\gr}{gr}
 
 \DeclareMathOperator{\Qcoh}{Qcoh}
\DeclareMathOperator{\coh}{coh} \DeclareMathOperator{\kr}{Ker}

 \DeclareMathOperator{\QGr}{QGr}
\DeclareMathOperator{\qgr}{qgr} \DeclareMathOperator{\Tors}{Tors}
\DeclareMathOperator{\tors}{tors} \DeclareMathOperator{\Tor}{Tor}
\DeclareMathOperator{\Proj}{\sf Proj}
\DeclareMathOperator{\perf}{\cc D_{per}}
 
\DeclareMathOperator{\Mod}{Mod} \DeclareMathOperator{\modd}{mod}

\newcommand {\lp}{\varinjlim}

\newcommand{\lra}[1]{\bl{#1}\longrightarrow\relax}
\newcommand{\bl}[1]{\buildrel #1\over}
\newcommand{\cc}{\mathcal}
\newcommand{\ps}{\oplus}

\newcommand{\homs}{\cc Hom}

\newcommand{\wt}{\widetilde}
\newcommand{\ifff}{if and only if }
\newcommand{\rfp}{\modd R}

\newcommand{\bb}{\mathbb}

\newcommand{\Rfp}{\Mod R}

\newtheorem{thm}{Theorem}[section]
\newtheorem{prop}[thm]{Proposition}
\newtheorem{cor}[thm]{Corollary}
\newtheorem{lem}[thm]{Lemma}

\newtheorem*{theo}{Theorem}

\newtheorem*{rem}{Remark}

\newtheorem*{sublem}{Sublemma}
\newtheorem*{defs}{Definition}

\numberwithin{equation}{section}

\makeatother

\begin{document}

\footskip30pt
%\baselineskip=1.2\baselineskip

\title{Reconstructing projective schemes from Serre subcategories}
\author{Grigory Garkusha}
\thanks{This paper was written during the visit of the first author to the
University of Manchester supported by the MODNET Research Training
Network in Model Theory. He would like to thank the University for
the kind hospitality.}
\address{Department of Mathematics, University of Wales Swansea, Singleton Park, SA2 8PP Swansea, UK}
\email{G.Garkusha@swansea.ac.uk}

\urladdr{homotopy.nm.ru}

\author{Mike Prest}
\address{School of Mathematics, University of Manchester, Oxford Road, M13~9PL Manchester, UK}
\date{August 23, 2006}
\email{mprest@maths.man.ac.uk}


\urladdr{www.maths.man.ac.uk/$\sim$mprest}

\begin{abstract}
Given a positively graded commutative coherent ring $A=\ps_{j\geq
0}A_j,$ finitely generated as an $A_0$-algebra, a bijection between
the tensor Serre subcategories of $\qgr A$ and the set of all
subsets $Y\subseteq\Proj A$ of the form $Y=\bigcup_{i\in\Omega}Y_i$
with quasi-compact open complement $\Proj A\setminus Y_i$ for all
$i\in\Omega$ is established. To construct this correspondence,
properties of the Ziegler and Zariski topo\-logies on the set of
isomorphism classes of indecomposable injective graded modules are
used in an essential way. Also, there is constructed an isomorphism
of ringed spaces
   $$(\Proj A,\cc O_{\Proj A})\lra{\sim}(\spec(\qgr A),\cc O_{\qgr A}),$$
where $(\spec(\qgr A),\cc O_{\qgr A})$ is a ringed space
associated to the lattice $L_{\serre}(\qgr A)$ of tensor Serre
subcategories of $\qgr A$.
\end{abstract}

%\dedicatory{}

\maketitle

\thispagestyle{empty} \pagestyle{plain}

\section{Introduction}

In his celebrated work on abelian categories P.~Gabriel~\cite{Ga}
proved that any noetherian scheme $X$ can be reconstructed uniquely
up to isomorphism from the category, $\Qcoh X$, of quasi-coherent
sheaves over $X$. This reconstruction result has been generalized to
arbitrary schemes by A.~Rosenberg in~\cite{R}.

Another result of Gabriel~\cite[VI.2, Prop.~4]{Ga} states that for
any noetherian scheme $X$ the assignments
   \begin{equation}\label{gabr}
    \coh X\supseteq\cc D\mapsto\bigcup_{x\in\cc D}\supp_X(x)\quad\text{and}\quad X\supseteq
    U\mapsto \{x\in\coh X\mid\supp_X(x)\subseteq U\}
   \end{equation}
induce bijections between
\begin{enumerate}
 \item the set of all Serre subcategories of $\coh X$, and
 \item the set of all subsets $U\subseteq X$ of the form
       $U=\bigcup_{i\in \Omega} Y_i$ where, for all $i \in \Omega$, $Y_i$
       has quasi-compact open complement
       $X\setminus Y_i$.
\end{enumerate}
As a consequence of this result, $X$ can be reconstructed from its
abelian category, $\coh X$, of coherent sheaves  (see
Buan-Krause-Solberg~\cite[Sec.~8]{BKS}).

Given a quasi-compact, quasi-separated scheme $X$, let $\perf(X)$
denote the derived category of perfect complexes. It comes equipped
with a tensor product $\otimes:=\otimes^L_{\cc O_X}$. A thick
triangulated subcategory $\cc T$ of $\perf(X)$ is said to be a
tensor subcategory if for every $E \in \perf(X)$ and every object
$A\in \cc T$, the tensor product $E\otimes A$ also is in $\cc T$.
Thomason~\cite{T} establishes a classification similar
to~\eqref{gabr} for tensor thick subcategories of $\perf(X)$ in
terms of the topology of $X$. Hopkins and Neeman (see~\cite{Hop,N})
did the case where $X$ is affine and noetherian.

Based on Thomason's classification theorem, Balmer~\cite{B1}
reconstructs the noetherian scheme $X$ from the tensor thick
triangulated subcategories of $\perf(X)$. This result has been
generalized to quasi-compact, quasi-separated schemes by
Buan-Krause-Solberg~\cite{BKS}.

In his fundamental paper~\cite{Se}, Serre proved a theorem which
describes the quasi-coherent sheaves on a projective scheme in terms
of graded modules as follows. Let $A$ be a finitely generated graded
algebra over a field $k$, and let $X=\Proj A$ be the associated
projective scheme. Let $\coh X$ denote the category of coherent
sheaves on $X$, and let $\cc O_X(n)$ denote the $n$th power of the
twisting sheaf on $X$. Define a functor $\Gamma_*:\coh X\to\qgr A$,
where $\qgr A$ is the category of finitely presented graded
$A$-modules modulo the shifts of $A_0$, by
   $$\Gamma_*(\cc F)=\bigoplus_{d=-\infty}^\infty H^0(X,\cc F\otimes\cc
   O_X(d)).$$
Serre's theorem asserts that if $A$ is generated over $k$ by
finitely many elements of degree 1 then $\Gamma_*$ defines an
equivalence of categories $\coh X\to\qgr A$.

Let $A$ be a coherent graded commutative ring $A=\ps_{j\geq 0}A_j$
which is finitely generated as an $A_0$-algebra. The purpose of this
paper is to establish a classification similar to~\eqref{gabr} for
tensor Serre subcategories of $\qgr A$ in terms of $\Proj A$. More
precisely, we demonstrate the following.

\begin{theo}[Classification]\label{tol}%\marginpar{tol}
Let $A$ be a coherent graded commutative ring which is finitely
generated as an $A_0$-algebra. The assignments
   $$\qgr A\supseteq\cc S\mapsto\bigcup_{M\in\cc S}\supp_A(M)\textrm{ and }
     \Proj A\supseteq U\mapsto\{M\in\qgr A\mid\supp_A(M)\subseteq U\}$$
induce bijections between
\begin{enumerate}
 \item the set of all tensor Serre subcategories of $\qgr A$, and
 \item the set of all subsets $U\subseteq\Proj A$ of the form
       $U=\bigcup_{i\in \Omega} Y_i$ with quasi-compact open complement
       $X\setminus Y_i$ for all $i\in \Omega$,
\end{enumerate}
where $\supp_A(M)=\{P\in\Proj A\mid M_P\neq 0\}$.
\end{theo}

Following Buan-Krause-Solberg~\cite{BKS} we consider the lattice
$L_{\serre}(\qgr A)$ of tensor Serre subcategories of $\qgr A$ and
its prime ideal spectrum $\spec(\qgr A)$. This space comes naturally
equipped with a sheaf of rings $\cc O_{\qgr A}$. The following
result says that the scheme $(\Proj A,\cc O_{\Proj A})$ is
isomorphic to $(\spec(\qgr A),\cc O_{\qgr A})$.

\begin{theo}[Reconstruction]\label{coh}
Let $A$ be a coherent graded ring which is finitely generated as an
$A_0$-algebra. There is a natural isomorphism
   $$f:(\Proj A,\cc O_{\Proj A})\lra{\sim}(\spec(\qgr A),\cc O_{\qgr A})$$
of ringed spaces.
\end{theo}

This theorem says that the abelian category $\qgr A$ contains all
the necessary information to reconstruct the projective scheme
$(\Proj A,\cc O_{\Proj A})$.

Our approach, similar to that used in \cite{GP}, makes use of
results on the Ziegler topology and its dual which had their origins
in the model theory of modules \cite{Z,Pr2}.

Throughout this paper we fix a positively graded commutative ring
$A=\ps_{j\geq 0}A_j$ with unit.

\section{Preliminaries}

In this section we recall some basic facts about graded rings and
modules.

\begin{defs}{\rm
A {\it (positively) graded ring\/} is a ring $A$ together with a
direct sum decomposition $A=A_0\ps A_1\ps A_2\ps\cdots$ as abelian
groups, such that $A_iA_j\subset A_{i+j}$ for $i,j\geq 0$. A {\it
homogeneous element\/} of $A$ is simply an element of one of the
groups $A_j$, and a {\it homogeneous ideal\/} of $A$ is an ideal
that is generated by homogeneous elements. A {\it graded
$A$-module\/} is an $A$-module $M$ together with a direct sum
decomposition $M=\ps_{j\in\bb Z}M_j$ as abelian groups, such that
$A_iM_j\subset M_{i+j}$ for $i\geq 0,j\in\bb Z$. One calls $M_j$ the
$j$th {\it homogeneous component of $M$}. The elements $x\in M_j$
are called {\it homogeneous (of degree $j$)}.

Note that $A_0$ is a commutative ring with $1\in A_0$, that all
summands $M_j$ are $A_0$-modules, and that $M=\ps_{j\in\bb Z}M_j$ is
a direct sum decomposition of $M$ as an $A_0$-module.

Let $A$ be a graded ring. The {\it category of graded $A$-modules},
denoted by $\Gr A$, has as objects the graded $A$-modules. A {\it
morphism\/} of graded $A$-modules $f:M\to N$ is an $A$-module
homomorphism satisfying $f(M_j)\subset N_j$ for all $j\in\bb Z$. An
$A$-module homomorphism which is a morphism in $\Gr A$ will be
called {\it homogeneous}.

Let $M$ be a graded $A$-module and let $N$ be a submodule of $M$.
$N$ is called a {\it graded submodule\/} if it is a graded module
such that the inclusion map is a morphism in $\Gr A$. The graded
submodules of $A$ are called {\it graded ideals}. If $d$ is an
integer the {\it tail\/} $M_{\geq d}$ is the graded submodule of $M$
having the same homogeneous components $(M_{\geq d})_j$ as $M$ in
degrees $j\geq d$ and zero for $j<d$. We also denote the ideal
$A_{\geq 1}$ by $A_+$.

}\end{defs}

For $n\in\bb Z$, $\Gr A$ comes equipped with a shift functor
$M\mapsto M(n)$ where $M(n)$ is defined by $M(n)_j=M_{n+j}$. It is a
Grothendieck category with the generating family $\{A(n)\}_{n\in\bb
Z}$. The tensor product for the category of all $A$-modules induces
a tensor product on $\Gr A$: given two graded $A$-modules $M,N$ and
homogeneous elements $x\in M_i,y\in N_j$, set $\deg(x\otimes
y):=i+j$. We define the {\it homomorphism $A$-module\/} $\cc
Hom_A(M,N)$ as follows. In dimension $n\in\bb Z$, the group $\cc
Hom_A(M,N)_n$ is the group of graded $A$-module homomorphisms of
degree $n$, i.e.,
   $$\cc Hom_A(M,N)_n=\Gr A(M,N(n)).$$

We refer to a graded $A$-module $M$ as {\it finitely generated\/} if
it is a quotient of a free graded module of finite rank
$\bigoplus_{s=1}^nA(d_s)$ where $d_1,\ldots,d_s\in\bb Z$. $M$ is
{\it finitely presented\/} if there is an exact sequence
   $$\bigoplus_{t=1}^mA(c_t)\to\bigoplus
   _{s=1}^nA(d_s)\to M\to 0.$$
The full subcategory of graded finitely presented modules will be
denoted by $\gr A$. Note that any graded $A$-module is a direct
limit of finitely presented graded $A$-modules.

The graded ring $A$ is said to be {\it coherent\/} if every finitely
generated graded ideal of $A$ is finitely presented. It is easy to
see that $A$ is coherent \ifff the category $\gr A$ is abelian.

Let $E$ be any indecomposable injective graded $A$-module (we remind
the reader that the corresponding ungraded module, $\bigoplus_nE_n$,
need not be injective in the category of ungraded $A$-modules). Set
$P= P(E)$ to be the sum of annihilator ideals $\ann_A(x)$ of
non-zero homogeneous elements $x\in E$. Observe that each ideal
$\ann_A(x)$ is homogeneous. Since $E$ is uniform the set of
annihilator ideals of non-zero homogeneous elements of $E$ is closed
under finite sum so the only issue is whether the sum, $P(E)$, of
them all is itself one of these annihilator ideals.

Given a prime homogeneous ideal $P$, we use the notation $E_P$ to
denote the injective hull, $E(A/P)$, of $A/P$. Notice that $E_P$ is
indecomposable. We also denote the set of isomorphism classes of
indecomposable injective graded $A$-modules by $\inj A$.

\begin{lem}\label{9.2}%\marginpar{9.2}
If $E\in \inj A$ then $P(E)$ is a homogeneous prime ideal. If the
module $E$ has the form $E_P(n)$ for some prime homogeneous ideal
$P$ and integer $n$, then $P=P(E)$.
\end{lem}
\begin{proof}
The proof is similar to that of~\cite[9.2]{Pr2}.
\end{proof}

It follows from the preceding lemma that the map
   $$P\subset A\mapsto E_P\in\inj A$$
from the set of homogeneous prime ideals to $\inj A$ is injective.

\section{The Zariski and Ziegler topologies}

Given any graded ring $A$ and any homogeneous ideal $I$ of $A$, let
us set $D^{\rm m}(I)=\{E\in\inj A\mid \cc Hom_A(A/I,E)=0\}$ (``m''
for ``morphism''). Since $D^{\rm m}(I)\cap D^{\rm m}(J)=D^{\rm
m}(I\cap J)$ (for the non-immediate inclusion, note that any
morphism from $A/(I\cap J)$ to $E$ extends, by injectivity of $E$,
to one from $A/I\oplus A/J$) these form a basis for a topology on
$\inj A$.

\begin{lem}\label{unbasop}%\marginpar{unbasop}
If $A$ is any graded ring and $I$ is a finitely generated
homogeneous ideal of $A$ and $I=\sum _1^nI_i$, then $D^{\rm
m}(I)=\bigcup _1^n D^{\rm m}(I_i)$.
\end{lem}

\begin{proof} Suppose that $ E\notin \bigcup _1^nD^{\rm m}(I_i). $ Then for each
$i$ there is a non-zero morphism $ f_i:R/I_i(d_j)\to E. $ The
intersection of the images of these morphisms is non-zero and, since
$A$ is commutative, any element in this intersection is annihilated
by each $I_i$, hence by $I$, that is, $\cc Hom_A(R/I,E)\neq 0$ so
$E\notin D^{\rm m}(I)$, as required.
\end{proof}

There is another topology on $\inj A$. The collection of subsets
   $$[M]=\{E\in\inj A\mid \cc Hom_A(M,E)=0\}$$
with $M\in\gr A$ forms a basis of open subsets for the {\it Zariski
topology\/} on $\inj A$. This topological space will be denoted by
$\injzar A$. Observe that $[M]=[M(d)]$ for any integer $d$.

For a coherent graded ring $A$ the sets $D^{\rm m}(I)$ with $I$
running over finitely generated homogeneous ideals form a basis for
topology on $\inj A$ which we call the {\it fg-ideals topology}. We
use the fact that the intersection, $I\cap J$, of two finitely
generated ideals $I,J$ is finitely generated in a coherent graded
ring. Indeed, since $I\cap J=\kr(A\to A/I\ps A/J)$ and $\gr A$ is
abelian, then $I\cap J$ is in $\gr A$.

By definition, $D^{\rm m}(I)=[A/I]$ for $I$ a finitely generated
homogeneous ideal. Let $M$ be a finitely presented graded
$A$-module. It is finitely generated by $b_1,...,b_n$ say,
$\deg(b_j)=d_j$. Set $M_k=\sum _{j\leq k}b_jA$, $M_0=0$. Each factor
$C_j=M_j/M_{j-1}$ is cyclic and, we claim, $[M]=[C_1]\cap ...\cap
[C_n]$. For, if there is a non-zero morphism from $C_j$ to $E(n)$,
$n\in\bb Z$, then, by injectivity of $E(n)$, this extends to a
morphism from $M/M_{j-1}$ to $E(n)$ and hence there is induced a
non-zero morphism from $M$ to $E(n)$. Conversely, if $f:M\to E(n)$
is non-zero let $j$ be minimal such that the restriction of $f$ to
$M_j$ is non-zero. Then $f$ induces a non-zero morphism from $C_j$
to $E(n)$. Since each $C_j$ is cyclic and finitely presented there
are finitely generated ideals $I_j$, $1\leq j\leq n$, such that
$C_j\cong A/I_j(d_j)$. It follows that each $[C_j]$ coincides with
$D^{\rm m}(I_j)$, and hence $[M]=D^{\rm m}(I)$ with
$I=\bigcap_{1\leq j\leq n}I_j$ finitely generated. Thus we have
shown the following

\begin{prop}\label{aaa}%\marginpar{aaa}
Given a coherent graded ring $A$, the Zariski topology on $\inj A$
coincides with the fg-ideals topology.
\end{prop}

Let us consider the collection of sets
   $$(M)=\injzar A\setminus[M]=\{E\in\inj A\mid \cc Hom_A(M,E)\neq 0\}$$
with $M\in\gr A$ and set
   $$O(M)=\{E\in\inj A\mid\Gr A(M,E)\neq 0\}.$$
Obviously, $(M)=\bigcup_{n\in\bb Z}O(M(n))$.

\begin{prop}\label{he}%\marginpar{he}
Let $A$ be a coherent graded ring. The collection of sets $(M)$ with
$M\in\gr A$ forms a basis of quasi-compact open sets for a topology,
called the {\it Ziegler topology\/}, on $\inj A$. This topological
space will be denoted by $\injzg A$.
\end{prop}

\begin{proof}
Since $A$ is coherent by assumption, so $\Gr A$ is a locally
coherent category, the collection of sets $O(M)$ with $M\in\gr A$
forms a basis of quasi-compact open sets for a topology on $\inj A$
(see~\cite{H,Kr1}). Given $L,M\in\gr A$, it follows that $O(L)\cap
O(M)=\bigcup_\alpha O(K_\alpha)$ for some $K_\alpha\in\gr A$. We
have
   \begin{gather*}(L)
    \cap(M)=[\bigcup_{n\in\bb Z}O(L(n))]\cap[\bigcup_{n\in\bb Z}O(M(n))]=\\
     =\bigcup_{n\in\bb Z}[O(L(n))\cap O(M(n))] \text{ which clearly equals } \bigcup_{n,\alpha} O(K_\alpha(n))
     =\bigcup_{\alpha}(K_\alpha).
   \end{gather*}
Since each $O(M)$ is quasi-compact and the translates of any cover
of $O(M)$ cover $(M)$, then so is $(M)$ for the Ziegler topology on
$\inj A$.
\end{proof}

If $\cc A$ is an abelian category then a Serre subcategory is a full
subcategory $\cc S$ such that if $0\to A\to B\to C\to 0$ is a short
exact sequence in $\cc A$ then $B\in \cc S$ \ifff $A, C \in \cc S$.
Given a subcategory $\cc X$ in $\gr A$ with $A$ graded coherent, we
may consider the smallest Serre subcategory of $\gr A$ containing
$\cc X$. This Serre subcategory we denote, following Herzog
\cite{H}, by
   $$\surd\cc X=\bigcap\{\cc S\subseteq\gr A\mid\cc S\supseteq\cc X\textrm{ is Serre}\}.$$
There is an explicit description of $\surd\cc X$.

\begin{prop}\cite[3.1]{H} \label{bbb}%\marginpar{bbb}
Let $A$ be a graded coherent ring and let $\cc X$ be a subcategory
of $\gr A$. A graded finitely presented module $M$ is in $\surd\cc
X$ \ifff there is a finite filtration of $M$ by graded finitely
presented submodules
   $$M=M_0\geq M_1\geq\cdots\geq M_n=0$$
and, for each $i<n$, there is $N_i\in\cc X$ and there are graded
finitely presented submodules
   $$N_i\geq N_{i_1}\geq N_{i_2}$$ such that $M_i/M_{i+1}\cong N_{i_1}/N_{i_2}$.
\end{prop}

Given a subcategory $\cc X$ of $\gr A$ denote by
   $$[\cc X]=\{E\in\inj A\mid\cc Hom_A(M,E)=0\textrm{ for all $M\in\cc X$}\}.$$
We shall also write $(\cc X)$ to denote $\inj A\setminus[\cc X]$.

\begin{cor}\label{ccc}\cite[3.3]{H} %\marginpar{ccc}
Given a graded coherent ring $A$ and $\cc X\subseteq\gr A$, we have
$[\cc X]=[\surd\cc X]$ and $(\cc X)=(\surd\cc X)$.
\end{cor}

\begin{proof}
This immediately follows from Proposition~\ref{bbb} and the fact
that the functor $\Gr A(-,E)$ with $E$ graded injective preserves
exact sequences.
\end{proof}

Let $A$ be a graded coherent ring. A {\it tensor Serre
subcategory\/} of $\gr A$ (or $\Gr A$) is a Serre subcategory $\cc
S\subset\gr A$ (or $\Gr A$) such that for any $X\in\cc S$ and any
$Y\in\gr A$ the tensor product $X\otimes Y$ is in $\cc S$.

\begin{lem}\label{ccc1}%\marginpar{ccc1}
Let $A$ be a graded coherent ring. $\cc S$ is a tensor Serre
subcategory of $\gr A$ \ifff it is closed under shifts of objects,
i.e. $X\in\cc S$ implies $X(n)\in\cc S$ for any $n\in\bb Z$.
\end{lem}

\begin{proof}
Suppose that $\cc S$ is a tensor Serre subcategory of $\gr A$. Then
it is closed under shifts of objects, because $X(n)\cong X\otimes
A(n)$.

Assume the converse. Let $X\in\cc S$ and $Y\in\gr A$. Then there is
a surjection $\ps_1^nA(n_i)\bl f\twoheadrightarrow Y$. It follows
that $1_X\otimes f:\ps_1^nX(n_i)\to X\otimes Y$ is a surjection.
Since each $X(n_i)$ belongs to $\cc S$ then so does $X\otimes Y$.
\end{proof}

\begin{prop}\label{bbb1}%\marginpar{bbb1}
Let $A$ be a graded coherent ring. The maps
   $$U\bl\phi\mapsto\cc S_{U}=\{M\in\gr A\mid (M)\subset U\}$$
and
   $$\cc S\bl\psi\mapsto U_{\cc S}=\bigcup_{M\in\cc S}(M)$$
induce a 1-1 correspondence between the lattices of open sets of
$\injzg A$ and tensor Serre subcategories of $\gr A$.
\end{prop}

\begin{proof}
By~\cite[3.8]{H} and~\cite[4.2]{Kr1} the maps
   $$\gr A\supseteq\cc S\mapsto O=\bigcup_{M\in\cc S}O(M)$$
and
   $$O\mapsto\cc S_{O}=\{M\in\gr A\mid O(M)\subset O\}$$
induce a 1-1 correspondence between the Serre subcategories of $\gr
A$ and the open sets $O$ of $\inj A$ for the topology defined by the
sets $O(M), M\in\gr A$ (see above). Our assertion now follows from
the fact that $(M)=\cup_{n\in\bb Z}O(M(n))$ and Lemma~\ref{ccc1}.
\end{proof}

Recall that for any homogeneous ideal $I$ of a graded ring, $A$, and
homogeneous $r\in A$ we have an isomorphism
$A/(I:r)\cong(rA+I)/I(\deg(r))$, where $(I:r)$ is the homogeneous
ideal $\{s\in A\mid rs\in I\}$, induced by sending $1+(I:r)$ to
$r+I$.

The next result is the ``nerve" in our analysis.

\begin{thm}\textrm{\em(cf. Prest~\cite[9.6]{Pr2})} \label{ddd}%\marginpar{ddd}
Let $A$ be a commutative coherent graded ring, let $E$ be an
indecomposable injective graded $A$-module and let $P(E)$ be the
prime ideal defined before. Then $E$ and $E_{P(E)}$ are
topologically indistinguishable in $\injzg R$ and hence also in
$\injzar R$.
\end{thm}

\begin{proof} Let $I$ be such that $E=E(A/I)(d)$. For each homogeneous $r\in(A\setminus
I)_n$ we have, by the remark just above, that the annihilator of
$r+I\in E$ is $(I:r)$ and so, by definition of $P(E)$, we have
$(I:r)\leq P(E)$. The natural projection $(rA+I)/I(d)\cong
(A/(I:r))(d-n)\to (A/P(E))(d-n)$ extends to a morphism from $E$ to
$E_{P(E)}(d-n)$ which is non-zero on $r+I$. Forming the product of
these morphisms as $r$ varies over homogeneous elements in
$A\setminus I$, we obtain a morphism from $E$ to a product of
appropriately shifted copies of $E_{P(E)}$ which is monic on $A/I$
and hence is monic. Therefore $E$ is a direct summand of a product
of appropriately shifted copies of $E_{P(E)}$ and so $E\in(M)$
implies $E_{P(E)}\in(M)$, where $M\in\gr A$.

For the converse, take a basic Ziegler-open neighbourhood of
$E_{P(E)}$. This has the form $(M)$ for a finitely presented graded
module $M$. Now, $E_{P(E)}\in (M)$ means that there is a non-zero
morphism $f:M(d)\longrightarrow E_{P(E)}$ for some integer $d$. We
can construct a pullback diagram
   $$\xymatrix{L\ar[r]\ar@{ >->}[d]&A/P(E)\ar@{ >->}[d]\\
               M(d)\ar[r]^f&E_{P(E)}}$$
in which the vertical arrows are monic. Let $0\ne K\subset L$ be a
finitely generated submodule of $L$ such that the restriction,
$f'$, of $f$ to $K$ is non-zero. Notice that the image of $f'$ is
contained in $A/P(E)$. $K$ is in $\gr A$ because it is a finitely
generated submodule of $M(d)\in\gr A$ and $A$ is coherent by
assumption. Since $A/P(E)=\varinjlim A/I_\lambda$, where
$I_\lambda$ ranges over the annihilators of non-zero homogeneous
elements of $E$, and $K$ is finitely presented, $f'$ factorises
through one of the maps $A/I_{\lambda_0}\to A/P(E)$
($I_{\lambda_0}=\ann_A(x)$ for some homogeneous $x\in E$). In
particular, there is a non-zero morphism $K\to E(\deg(x))$ and
hence, by injectivity of $E$, an extension to a morphism
$M(d-\deg(x))\to E$, showing that $E\in (M)$, as required.
\end{proof}

\section{Torsion modules and the category $\QGr A$}

Throughout this section the graded ring $A$ is supposed to be
coherent and the homogeneous ideal $A_+\subset A$ is supposed to be
finitely generated. This is equivalent to saying that $A$ is
coherent and a finitely generated $A_0$-algebra. In this section we
introduce the category $\QGr A$ ($\qgr A$) which is analogous to the
category of quasi-coherent (coherent) sheaves on a projective
variety. The non-commutative analog of the category $\QGr A$ plays a
prominent role in ``non-commutative projective geometry" (see, e.g.,
\cite{AZ,S,V}). We use the general theory of locally coherent
Grothendieck categories and their localizations (see~\cite{H,Kr1}
for details).

Recall that a Serre subcategory $\cc S$ of $\Gr A$ is {\it
localizing\/} if it is closed under taking direct limits. A
localizing subcategory $\cc S$ is said to be {\it of finite type\/}
if the canonical functor from the quotient category $\Gr A/\cc
S\to\Gr A$ respects direct limits, equivalently if $\cc S$ is the
closure under direct limits of a Serre subcategory of $\gr A$.

Let us consider the commutative ring $A_0$ as a graded $A$-module.
It is isomorphic to $A/A_+$ and belongs to $\gr A$, because $A_+$ is
finitely generated by assumption. Put $\tors A=\surd
\{A_0(d)\}_{d\in\bb Z}$, the tensor Serre subcategory in $\gr A$
generated by the shifts of $A_0$, and set $\Tors A=\{\lp T_i\mid
T_i\in\tors A\}$. Then $\Tors A$ is a localizing subcategory of
finite type in $\Gr A$ and $\tors A=\Tors A\cap\gr A$
(see~\cite[2.8]{H}, \cite[2.8]{Kr1}). We refer to the objects of
$\Tors A$ as {\it torsion graded modules}.

\begin{lem}\label{dd1}%\marginpar{dd1}
Under the assumptions above the graded module $A/A_{\geq n}(d)$ is
torsion for any $n\geq 1$ and $d\in\bb Z$.
\end{lem}

\begin{proof}
The graded module $A/A_{\geq 1}(d)$ is torsion by assumption. We
proceed by induction on $n$. Suppose $A/A_{\geq n}(d)$ is torsion.
We want to check that $A/A_{\geq n+1}(d)$ is torsion.

Consider the short exact sequence in $\Gr A$
   $$0\to A_{\geq n}/A_{\geq n+1}\to A/A_{\geq n+1}\to A/A_{\geq n}\to 0.$$
Clearly $A_{\geq n}/A_{\geq n+1}$ is an epimorphic image of the
torsion graded module $\ps_{x\in A_{n}}A_0(n)$, hence is torsion
itself. Since $A_{\geq n}/A_{\geq n+1}$ is torsion and $\Tors A$ is
closed under extensions, we see that $A/A_{\geq n+1}$ is torsion as
well and therefore so is each $A/A_{\geq n+1}(d)$, $d\in\bb Z$.
\end{proof}

Let $\QGr A=\Gr A/\Tors A$. We define $\tau$ as the functor which
assigns to a graded $A$-module its maximal torsion module. By
$\pi:\Gr A\to\QGr A$ we denote the quotient functor. By standard
localization theory $\pi$ is exact and respects direct limits. We
denote the fully faithful right adjoint to $\pi$ by $\omega$ and we
denote the composition $\omega\pi$ by $Q$. Since $\pi\omega$ is the
identity, it follows that $Q^2=Q$. A graded module is said to be
{\it $\Tors$-closed\/} if it has the form $Q(M)$ for some $M\in\Gr
A$. We shall identify $\QGr A$ with the full subcategory of
$\Tors$-closed modules. For any $X\in\Gr A$ and $Y\in\QGr A$ there
is an isomorphism
   $$\Gr A(X,Y)\cong\QGr A(Q(X),Y)$$
natural both in $X$ and $Y$. The shift functor $M\mapsto M(n)$
defines a shift functor on $\QGr A$ for which we shall use the same
notation. Observe that the functors $\tau, Q$ commute with the shift
functor. Finally we shall write $\cc O=Q(A)$.

Put $\qgr A=\gr A/\tors A$. It is easy to see that the obvious
functor $\qgr A\to\QGr A$ is fully faithful. We shall identify $\qgr
A$ with its image in $\QGr A$. $\qgr A$ is an abelian category and
equals $\surd\{\cc O(d)\}_{d\in\bb Z}$, the smallest Serre
subcategory of $\qgr A$ containing $\{\cc O(d)\}_{d\in\bb Z}$.
Moreover, every object of $\QGr A$ can be written as a direct limit
$\lp M_i$ of objects $M_i\in\qgr A$. It follows from~\cite[2.16]{H}
that $M\in\qgr A$ \ifff it has the form $Q(L)$ for some $L\in\gr A$.
Therefore every $M\in\qgr A$ has a presentation in $\qgr A$
   $$\bigoplus_{t=1}^m\cc O(c_t)\to\bigoplus_{s=1}^n\cc O(d_s)\to M\to 0.$$

The tensor product in $\Gr A$ induces a tensor product in $\QGr A$,
denoted by $\boxtimes$. More precisely, one sets
   $$X\boxtimes Y:=Q(X\otimes Y)$$
for any $X,Y\in\QGr A$.

\begin{lem}\label{d2}%\marginpar{d2}
Given $X,Y\in\Gr A$ there is a natural isomorphism in $\QGr A$:
$Q(X)\boxtimes Q(Y)\cong Q(X\otimes Y)$. Moreover, the functor
$-\boxtimes Y:\QGr A\to\QGr A$ is right exact and preserves direct
limits.
\end{lem}

\begin{proof}
By standard localization theory (see~\cite{Ga}) there is a long
exact sequence
   $$0\to T\to X\lra{\lambda_X}Q(X)\to T'\to 0,$$
where $T=\tau(X)$, the largest torsion submodule of $X$, $T'\in\Tors
A$, and $Q(X)$ is the maximal essential extension of $X':=X/T$ such
that $Q(X)/X'\in\Tors A$. We have an exact sequence in $\Gr A$:
   $$T\otimes Y\to X\otimes Y\bl\ell\to X'\otimes Y\to 0.$$
Since $\Tors$ is a tensor Serre subcategory, we have $T\otimes
Y\in\Tors$. Therefore $Q(\ell)$ is an isomorphism.

On the other hand, one has an exact sequence
   \begin{equation}\label{kjh}
    \cdots\to\Tor_1(T',Y)\to X'\otimes Y\bl\iota\to Q(X)\otimes Y\to T'\otimes Y\to 0
   \end{equation}
with $T'\otimes Y\in\Tors A$.

We claim that $\Tor_1(T',Y)\in\Tors A$. Indeed, choose a free
resolution $F_*$
   $$\cdots\to F_1\to F_0\to Y\to 0$$
of $Y$. Then the homology groups of the complex $T'\otimes F_*$ are
$\Tor_{n\geq 0}(T',Y)$. But $T'\otimes F_*$ is a complex in $\Tors
A$, hence its homology belongs to $\Tors A$.

Applying the exact functor $Q$ to sequence~\eqref{kjh} we infer that
$Q(\iota)$ is an isomorphism. It follows that $Q(\lambda_X\otimes
1_Y)=Q(\iota)Q(\ell)$ is an isomorphism as well. We have:
   $$Q(X\otimes Y)\cong Q(Q(X)\otimes Y)$$
which, by the same reasoning, is isomorphic to $$Q(Q(X)\otimes
Q(Y))=Q(X)\boxtimes Q(Y).$$ We leave the reader to check that the
functor $-\boxtimes Y$ is right exact and preserves direct limits.
\end{proof}

As a consequence of the preceding lemma we get an isomorphism
$X(d)\cong\cc O(d)\boxtimes X$ for any $X\in\QGr A$ and $d\in\bb Z$.

The notion of a tensor Serre subcategory of $\qgr A$ (with respect
to the tensor product $\boxtimes$) is defined similar to tensor
Serre subcategories of $\gr A$. The proof of the next lemma is like
that of Lemma~\ref{ccc1} (also use Lemma~\ref{d2}).

\begin{lem}\label{ccc2}%\marginpar{ccc2}
$\cc S$ is a tensor Serre subcategory of $\qgr A$ \ifff it is closed
under shifts of objects, i.e. $X\in\cc S$ implies $X(n)\in\cc S$ for
any $n\in\bb Z$.
\end{lem}

\section{Some properties of $\Proj A$}

Recall that the projective scheme $\Proj A$ is a topological space
whose points are the graded prime ideals not containing $A_+$. The
topology of $\Proj A$ is defined by taking the closed sets to be the
sets of the form $V(I)=\{P\in\Proj A\mid P\supseteq I\}$ for some
homogeneous ideal $I$ of $A$. We set $D(I):=\Proj A\setminus V(I)$.

We also recall from~\cite{Hoc} that a topological space is {\it
spectral\/} if it is $T_0$ and quasi-compact, the quasi-compact open
subsets are closed under finite intersections and form an open
basis, and every non-empty irreducible closed subset has a generic
point.

\begin{prop}\label{d3}%\marginpar{d3}
Let $A$ be a graded ring which is finitely generated as an
$A_0$-algebra. Then the space $\Proj A$ is spectral.
\end{prop}

\begin{proof}
Let $P,Q$ be two different graded prime ideals of $A$. Without loss
of generality we may assume that there is a homogeneous element
$a\in P$ such that $a\notin Q$. It follows that $P\notin D(a)$ but
$Q\in D(a)$. Therefore $\Proj A$ is $T_0$.

Below we shall need the following

\begin{sublem}
A graded ideal $P$ is prime \ifff for any two homogeneous elements
$x,y\in A$ the condition $xy\in P$ implies that $x\in P$ or $y\in
P$.
\end{sublem}

\begin{proof}
Assume that for any two homogeneous elements $x,y\in A$ the
condition $xy\in P$ implies $x\in P$ or $y\in P$. Let $a,b\in A$ be
such that $ab\in P$. We write $a=\sum_ia_i$, $a_i\in A_i$, and
$b=\sum_jb_j$, $b_j\in A_j$. Assume that $a\notin P$ and $b\notin
P$. Then there exist integers $p,q$ such that $a_p\notin P$, but
$a_i\in P$ for $i<p$, and $b_q\notin P$, but $b_j\in P$ for $j<q$.
The $(p+q)$th homogeneous component of $ab$ is
$\sum_{i+j=p+q}a_ib_j$. Thus $\sum_{i+j=p+q}a_ib_j\in P$, since $P$
is graded. All summands of this sum, except possibly $a_pb_q$,
belong to $P$, and so it follows that $a_pb_q\in P$ as well. We
conclude by assumption that $a_p\in P$ or $b_q\in P$, a
contradiction.
\end{proof}

Let $a$ be a non-zero homogeneous element of $A_+$. Let us show that
$D(a)$ is quasi-compact. For this, we consider a cover
$D(a)=\bigcup_\Lambda D(I_\lambda)$ of $D(a)$ by open sets
$D(I_\lambda)$. Assume $D(a)\ne D(I_{\lambda_1})\cup\cdots\cup
D(I_{\lambda_n})$ for any $\lambda_1,\ldots,\lambda_n\in\Lambda$.
Set $I:=\sum_\Lambda I_\lambda$. Then $a\notin I$ because otherwise
$a\in I_{\lambda_1}+\cdots+I_{\lambda_n}$ for some
$\lambda_1,\ldots,\lambda_n\in\Lambda$ and then
$D(a)=D(I_{\lambda_1})\cup\cdots\cup D(I_{\lambda_n})$. It also
follows that $a^t\notin I$ for any $t$ and $I\varsupsetneq A_+$.

\begin{sublem}
Let $Q$ be a graded ideal with $a^t\notin Q$ for all $t$, $Q\geq I$,
and let $Q$ be maximal such. Then $Q$ is prime.
\end{sublem}

\begin{proof}
By the sublemma above it is enough to check that for any two
homogeneous elements $b,c\in A$ the condition $bc\in Q$ implies
$b\in Q$ or $c\in Q$. Let $b,c\in A$ be such that $bc\in Q$ and
$b,c\notin Q$. Then $a^t=q+br\in Q+bA$ and $a^s=q'+cr'\in Q+cA$ for
some $s,t\in\bb N$ and $q,q'\in Q$. We see that
$a^{t+s}=q''+bcr''\in Q$, a contradiction.
\end{proof}

Thus $a\notin Q$ hence $Q\in D(a)$ and so $Q\in D(I_{\lambda_0})$
for some $\lambda_0\in\Lambda$. It follows that $Q\varsupsetneq
I_{\lambda_0}$, a contradiction. So $D(a)$ is quasi-compact, and
hence every subset $D(J)$ with $J$ a finitely generated graded ideal
is quasi-compact and every quasi-compact open subset is of this
form.

Let $a_1,\ldots,a_n$ be generators for $A_+$. It follows that $\Proj
A = D(a_1)\cup\cdots\cup D(a_n)$ is quasi-compact. Also the
quasi-compact open subsets form an open basis since the set of them
is closed under finite intersections: given any finitely generated
graded ideals $I,J$ with generators $b_1,\ldots,b_l$ and
$c_1,\ldots,c_m$ respectively, we have $D(I)\cap
D(J)=\bigcup_{i,j}D(b_ic_j)$ and each $D(b_ic_j)$ is quasi-compact.

It remains to verify that every non-empty irreducible closed subset
$V$ has a generic point.

There exists a graded ideal $I$ such that $V=V(I)$. Without loss of
generality we may assume that $I=\surd I$, where $\surd I$ denotes
the graded ideal generated by the homogeneous elements $b$ such that
$b^t\in I$ for some $t$. If $I$ is prime then it is a generic point
of $V$. If $I$ is not prime there are homogeneous $a,b\in A\setminus
I$ such that $ab\in I$ (see the first sublemma above). Then
$V=V(aA+I)\cup V(bA+I)$, hence $V=V(aA+I)$ for $V$ is irreducible.
Therefore every graded prime ideal $P$ containing $I$ must contain
$a$.

Let $Q$ be a graded ideal with $a^t\notin Q$ for any $t$, $Q\geq I$,
and let $Q$ be maximal such. The proof of the second sublemma above
shows that $Q$ is prime. It follows that $Q\in V\setminus
V(aA+I)=\emptyset$, contradiction.
\end{proof}

In the remainder of this section $A$ is coherent and a finitely
generated $A_0$-algebra.  We denote by $\QInj A$ the subset of $\Inj
A$ consisting of torsion-free indecomposable injectives.

\begin{lem}\label{d1}%\marginpar{d1}
If $A$ is coherent and a finitely generated $A_0$-algebra then
$E\in\QInj A$ \ifff the graded prime ideal $P(E)$ is in $\Proj A$.
\end{lem}

\begin{proof}
By assumption $A_0 \in \gr A$.  Let $E\in\QInj A$. Then $\cc
Hom(A_0,E)=0$ and hence $\cc Hom(A_0,E_{P(E)})=0$ by
Theorem~\ref{ddd}. Therefore $P(E)\varsupsetneq A_+$ because
otherwise the composite map $A_0\cong A/A_+\twoheadrightarrow
A/P(E)\rightarrowtail E_{P(E)}$ would be non-zero. Thus
$P(E)\in\Proj A$.

Now assume the converse. Suppose that there is a non-zero map
$f:A_0\to E_{P(E)}(d)$ for some integer $d$. It follows that
$A_+\subset P(E)$, a contradiction. By Theorem~\ref{ddd} $\cc
Hom(A_0,E)=0$, hence $E$ has no torsion.
\end{proof}

\begin{cor}\label{d111}%\marginpar{d111}
For any power $A_+^t$ of $A_+$ the module $A/A_+^t(d), d\in\bb Z,$
is torsion.
\end{cor}

\begin{proof}
Let $E\in\QInj A$. Suppose $\cc Hom(A/A^t_+,E)\ne 0$. Since $A_+$ is
finitely generated, then so is $A^t_+$ and hence $A/A^t_+\in\qgr A$.
By Theorem~\ref{ddd} $\cc Hom(A/A^t_+,E_{P(E)})\ne 0$, and so
$P(E)\supset A_+^t$. It follows that $P(E)\supset A_+$ what is
impossible by the preceding lemma.
\end{proof}

It is useful to have the following characterization of torsion
modules.

\begin{prop}\label{d222}%\marginpar{d222}
A graded module $T$ is torsion \ifff every element $x\in T$ is
annihilated by some power $A_+^t$ of $A_+$, that is,
$\ann_A(x)\supseteq A_+^t$.
\end{prop}

\begin{proof}
Suppose first that every element $x\in T$ is annihilated by some
power $A_+^{t_x}$ of $A_+$. Then there is an epimorphism
   $$\bigoplus_{x\in T}A/A_+^{t_x}(\deg(x)) \to T.$$
The module on the left is torsion by Corollary~\ref{d111}, and
hence so is $T$.

Now observe that the full subcategory of $\QGr A$ consisting of
those modules whose elements are annihilated by powers of $A_+$ is
localizing. Moreover this subcategory contains $A/A_+$, hence all of
$\tors A$. Therefore it contains, hence, by the first paragraph,
coincides with, $\Tors A$.
\end{proof}

Given a prime ideal $P\in\Proj A$ and a graded module $M$, denote by
$M_P$ the homogeneous localization of $M$ at $P$. If $f$ is a
homogeneous element of $A$, by $M_f$ we denote the localization of
$M$ at the multiplicative set $S_f=\{f^n\}_{n\geq 0}$.

\begin{cor}\label{g111}%\marginpar{g111}
If $T$ is a torsion module then $T_P=0$ and $T_f=0$ for any
$P\in\Proj A$ and $f\in A_+$. As a consequence, $M_P\cong Q(M)_P$
and $M_f\cong Q(M)_f$ for any $M\in\Gr A$.
\end{cor}

\begin{proof}
Every homogeneous element of $T_P$ is of the form $x/s$, where $x\in
T,s\in A\setminus P$ are homogeneous elements. Since $P\in\Proj A$
there exist a homogeneous element $f\in (A\setminus P)\cap A_+$. By
Proposition~\ref{d222} $x$ is annihilated by some power $A_+^t$ of
$A_+$. It follows that $f^ts\cdot x/s=0$, and hence $T_P=0$. The
fact that $T_f=0$ is proved in a similar fashion.

Now let $X$ be a graded $A$-module. There is a long exact sequence
   $$0\to T\to X\lra{\lambda_X}Q(X)\to T'\to 0,$$
where $T,T'\in\Tors A$. We conclude that $(\lambda_X)_P$ and
$(\lambda_X)_f$ are isomorphisms.
\end{proof}

\section{The Zariski and Ziegler topologies on $\QInj A$}

Throughout this section the graded ring $A$ is supposed to be
coherent and a finitely generated $A_0$-algebra. The {\it Zariski
topology\/} on $\QInj A$ is given by the collection of sets
   $$[M]_q:=\{E\in\QInj A\mid\cc Hom_A(M,E)=0\}$$
with $M\in\qgr A$ forming a basis of open subsets. This topological
space will be denoted by $\qinjzar A$. Observe that $[M]_q=[M(d)]_q$
for any integer $d$. The Zariski topology on $\QInj A$ coincides
with the subspace topology induced by the Zariski topology on $\inj
A$, because
   $$\Gr A(L,E)\cong\QGr A(Q(L),E)$$
for any $L\in\gr A$ and hence $[L]\cap\QInj A=[Q(L)]_q$.

The Ziegler topology on $\QInj A$ is defined in a similar way.
Namely, set
   $$O_q(M):=\{E\in\QInj A\mid\QGr A(M,E)\ne 0\}$$
and
   $$(M)_q:=\bigcup_{d\in\bb Z}O_q(M(d))=\{E\in\QInj A\mid\cc Hom_A(M,E)\ne 0\}$$
with $M\in\qgr A$. The proof of Proposition~\ref{he} shows that the
collection of sets $(M)_q$ with $M\in\qgr A$ forms a basis of
quasi-compact open sets for the {\it Ziegler topology\/} on $\QInj
A$. This topological space will be denoted by $\qinjzg A$. Moreover,
the Ziegler topology on $\QInj A$ is compatible with the Ziegler
topology on $\inj A$, because $O(L)\cap\QInj A=O_q(Q(L))$ and
$(L)\cap\QInj A=(Q(L))_q$ for any $L\in\gr A$.

Given any homogeneous ideal $I$ of $A$, let us set $D^{\rm
m,q}(I)=\{E\in\QInj A\mid \cc Hom_A(A/I,E)=0\}$. The {\it
fg-topology\/} on $\QInj A$ is defined by the collection of (basic
open) subsets $D^{\rm m,q}(I)=[Q(A/I)]_q$ with $I$ a finitely
generated homogeneous ideal. We have $D^{\rm m}(I)\cap\QInj A=D^{\rm
m,q}(I)$. Moreover, the Zariski topology and the fg-ideals topology
for $\QInj A$ coincide.

Using Lemmas~\ref{9.2} and \ref{d1}, there is an embedding
   $$\alpha:\Proj A\to\QInj A,\ \ \ P\mapsto E_P.$$
We shall identify $\Proj A$ with its image in $\QInj A$.

\begin{thm}\label{vera}%\marginpar{vera}
The space $\Proj A$ is dense and a retract in $\qinjzar A$. A left
inverse to the embedding $\Proj A\hookrightarrow\qinjzar A$ takes an
indecomposable injective torsion-free graded module $E$ to the prime
ideal $P(E)$ which is the sum of annihilator ideals of non-zero
homogeneous elements of $E$. Moreover, $\qinjzar A$ is
quasi-compact, the basic open subsets $[M]_q, M\in\qgr A,$ are
quasi-compact, the intersection of two quasi-compact open subsets is
quasi-compact, and every non-empty irreducible closed subset has a
generic point.
\end{thm}

\begin{proof}
For any ideal $I$ of $A$ we have
   \begin{equation}\label{zzz}
    D^{\rm m,q}(I)\cap\Proj A=D(I).
   \end{equation}
From this relation, the fact that the Zariski topology and the
fg-ideals topology for $\QInj A$ coincide, and Theorem~\ref{ddd} it
follows that $\Proj A$ is dense in $\qinjzar A$ and that
$\alpha:\Proj A\to\qinjzar A$ is a continuous map.

It follows from Lemma~\ref{9.2} that
   $$\beta:\qinjzar A\to\Proj A,\ \ \ E\mapsto P(E),$$
is left inverse to $\alpha$. The relation~\eqref{zzz} implies that
$\beta$ is continuous as well. Thus $\Proj A$ is a retract of
$\qinjzar A$.

Let us show that each basic open set $[M]_q$, $M\in\qgr A$, is
quasi-compact (in particular $\qinjzar A=[0]_q$ is quasi-compact).
$[M]_q=D^{\rm m,q}(I)$ for some finitely generated graded ideal $I$
of $A$.

Let $D^{\rm m,q}(I)=\bigcup_{i\in\Omega}D^{\rm m,q}(I_i)$ with each
$I_i$ graded finitely generated. It follows from~\eqref{zzz} that
$D(I)=\bigcup_{i\in\Omega}D(I_i)$. Since $I$ is finitely generated,
$D(I)$ is quasi-compact in $\Proj A$ by the proof of
Proposition~\ref{d3}. We see that
$D(I)=\bigcup_{i\in\Omega_0}D(I_i)$ for some finite subset
$\Omega_0\subset\Omega$.

Assume $E\in D^{\rm m,q}(I)\setminus\bigcup_{i\in\Omega_0}D^{\rm
m,q}(I_i)$. It follows from Theorem~\ref{ddd} that $E_{P(E)}\in
D^{\rm m,q}(I)\setminus\bigcup_{i\in\Omega_0}D^{\rm m,q}(I_i)$. But
$E_{P(E)}\in D^{\rm m,q}(I)\cap\Proj
A=D(I)=\bigcup_{i\in\Omega_0}D(I_i)$, and hence it is in
$D(I_{i_0})=D^{\rm m,q}(I_{i_0})\cap\Proj A$ for some
$i_0\in\Omega_0$, a contradiction. So $[M]_q$ is quasi-compact.

It follows from above that every quasi-compact open subset in
$\qinjzar A$ is of the form $[M]_q$ with $M\in\qgr A$. Therefore the
intersection $[M]_q\cap[N]_q=[M\ps N]_q$ of two quasi-compact open
subsets is quasi-compact.

By Theorem~\ref{ddd}, the proof of Proposition~\ref{d3}, and
relation~\eqref{zzz} we see that a subset $V$ of $\qinjzar A$ is
Zariski-closed and irreducible \ifff there is a graded prime ideal
$Q$ of $A$ such that $V=\{E\mid P(E)\geq Q\}$. Theorem~\ref{ddd},
Lemma~\ref{9.2}, and \eqref{zzz} obviously imply that the point
$E_Q\in V$ is generic.
\end{proof}

\begin{rem}{\rm
A similar result for affine schemes $\spec R$ with $R$ a commutative
coherent ring has been shown in~\cite[Theorem~A]{GP}. }
\end{rem}

Given a spectral topological space, $X$, Hochster \cite{Hoc} endows
the underlying set with a new, ``dual", topology, denoted $X^*$, by
taking as open sets those of the form $Y=\bigcup_{i\in\Omega}Y_i$
where $Y_i$ has quasi-compact open complement $X\setminus Y_i$ for
all $i\in\Omega$. Then $X^*$ is spectral and $(X^*)^*=X$ (see
\cite[Prop.~8]{Hoc}). The spaces, $X$, which we consider here are
not in general spectral; nevertheless we make the same definition
and denote the space so obtained by $X^*$. We also write $\Proj^* A$
for $(\Proj A)^*$.

\begin{cor}\label{eee}%\marginpar{eee}
The following relations hold:
   $$\qinjzar A=(\qinjzg A)^* \textrm{ and }\qinjzg A=(\qinjzar A)^*.$$
\end{cor}

\begin{proof}
This follows from Theorem~\ref{vera} and the fact that a
Ziegler-open subset $ U$ is quasi-compact \ifff it is one of the
basic open subsets $(M)_q, M\in\qgr A$.
\end{proof}

The proof of the following statement literally repeats that of
Proposition~\ref{bbb1}.

\begin{prop}\label{bbb2}%\marginpar{bbb2}
The maps
   $$ U\bl\phi\mapsto\cc S_{ U}=\{M\in\qgr A\mid (M)_q\subset U\}$$
and
   $$\cc S\bl\psi\mapsto U_{\cc S}=\bigcup_{M\in\cc S}(M)_q$$
induce a 1-1 correspondence between the lattices of open sets of
$\qinjzg A$ and tensor Serre subcategories of $\qgr A$.
\end{prop}

\begin{lem}\label{lll}%\marginpar{lll}
The maps
   $$\Proj^* A\supseteq  U\bl\phi\mapsto\cc Q_{ U}=\{E\in\qinjzg A\mid P(E)\in U\}$$
and
   $$\qinjzg A\supseteq \cc Q\bl\psi\mapsto U_{\cc Q}=\{P(E)\in\Proj^* A\mid E\in\cc Q\}=
     \cc Q\cap\Proj^*A$$
induce a 1-1 correspondence between the lattices of open sets of
$\Proj^*A$ and those of $\qinjzg A$.
\end{lem}

\begin{proof}
First note that $E_P\in\cc Q_{ U}$ for any $P\in U$ (see
Lemmas~\ref{9.2} and~\ref{d1}). Let us check that $\cc Q_{ U}$ is
an open set in $\qinjzg A$. Given a graded ideal $I$ of $A$,
denote by $V(I):=\Proj A\setminus D(I)$ and $V^{\rm
m,q}(I):=\qinjzg A\setminus D^{\rm m,q}(I)$. By definition, each
$V(I)$ with $I$ a finitely generated graded ideal of $A$ is a
basic open set in $\Proj^*A$. It follows from~\eqref{zzz} that
   \begin{equation}\label{www}
    V(I)=V^{\rm m,q}(I)\cap\Proj^*A.
   \end{equation}
Every closed subset of $\Proj A$ with quasi-compact complement has
the form $V(I)$ for some finitely generated graded ideal, $I$, of
$A$ (see Proposition~\ref{d3}), so there are finitely generated
graded ideals $I_\lambda\subseteq A$ such that $ U=\bigcup_\lambda
V(I_\lambda)$. Since the points $E$ and $E_{P(E)}$ are,  by
Theorem~\ref{ddd}, indistinguishable in $\qinjzg A$ we see that $\cc
Q_{ U}=\bigcup_\lambda V^{\rm m,q}(I_\lambda)$, hence this set is
open in $\qinjzg A=(\qinjzar A)^*$.

The same arguments imply that $ U_{\cc Q}$ is open in $\Proj^*A$. It
is now easy to see that $ U_{\cc Q_{ U}}= U$ and $\cc Q_{U_{\cc
Q}}=\cc Q$.
\end{proof}

\section{Classifying Serre subcategories of $\qgr A$}

Throughout this section the graded ring $A$ is supposed to be
coherent and a finitely generated $A_0$-algebra.

If $M\in\qgr A$ we write
   $$\supp_A(M)=\{P\in\Proj^*A\mid M_P\neq 0\},$$
where $M_P$ stands for the homogeneous localization of the graded
module $M$ at graded prime $P$.

\begin{lem}\label{iii}%\marginpar{iii}
Let $M\in\qgr A$ and $E\in\QInj A$. Then $E\in(M)_q$ \ifff
$M_{P(E)}\ne 0$ (or equivalently $P(E)\in\supp_A(M)$).
\end{lem}

\begin{proof}
By Theorem~\ref{ddd}, $E\in(M)_q$ \ifff $E_{P(E)}\in(M)_q$. The
assertion now follows from the easy fact that $M_P=0$ \ifff $\cc
Hom_A(M,E_P)=0$.
\end{proof}

Given a subcategory $\cc X$ of $\qgr A$ denote by
   $$[\cc X]_q=\{E\in\QInj A\mid\cc Hom_A(M,E)=0\textrm{ for all $M\in\cc X$}\}.$$
We shall also write $(\cc X)_q$ to denote $\QInj A\setminus[\cc
X]_q$.

It follows from the preceding lemma that
   \begin{equation}\label{vvv}
    \supp_A(M)=(M)_q\cap\Proj^*A.
   \end{equation}
Hence $\supp_A(M)$ is an open set of $\Proj^*A$ by
Lemma~\ref{lll}. More generally we have for any $\cc
X\subseteq\qgr A$:
    $$\supp_A(\cc X):=\bigcup_{M\in\cc X}\supp_A(M)=(\cc X)_q\cap\Proj^*A.$$
The proof of Corollary~\ref{ccc} shows that $(\cc X)_q=(\surd\cc
X)_q$. Therefore,
   \begin{equation}\label{vv}
    \supp_A(\cc X)=\supp_A(\surd\cc X).
   \end{equation}

We shall write $L_{\open}(\Proj^*A)$, $L_{\serre}(\qgr A)$ to
denote:
\begin{itemize}
\item{$\diamond$} the lattice of all open subsets of $\Proj^*A$,

\item{$\diamond$} the lattice of all tensor Serre subcategories of $\qgr A$.
\end{itemize}

We are now in a position to demonstrate the main result of this
section classifying tensor Serre subcategories of $\qgr A$.

\begin{thm}[Classification]\label{tol}%\marginpar{tol}
The assignments
   $$\qgr A\supseteq\cc S\mapsto\bigcup_{M\in\cc S}\supp_A(M)\textrm{ and }
     \Proj A\supseteq U\mapsto\{M\in\qgr A\mid\supp_A(M)\subseteq U\}$$
induce a lattice isomorphism between $L_{\open}(\Proj^*A)$ and
$L_{\serre}(\qgr A)$.
\end{thm}

\begin{proof}
The map
   $$\tau:\qgr A\supseteq\cc S\mapsto\bigcup_{M\in\cc S}\supp_A(M)$$
factors as
   $$\qgr A\supseteq\cc S\bl\delta\mapsto\cc Q=\bigcup_{M\in\cc S}(M)_q
     \bl\psi\mapsto\bigcup_{M\in\cc S}\supp_A(M),$$
where $\psi$ is the map of Lemma~\ref{lll} (we have used here
relation~\eqref{vvv}). By Proposition~\ref{bbb2} the map $\delta$
induces a 1-1 correspondence between the tensor Serre subcategories
of $\qgr A$ and the open sets $\cc Q$ of $\qinjzg A$. By
Lemma~\ref{lll} the map $\psi$ induces a 1-1 correspondence between
the lattices of open sets of $\qinjzg A$ and those of $\Proj^*A$.
Therefore the map $\tau$ induces the desired 1-1 correspondence
between the tensor Serre subcategories of $\qgr A$ and the open sets
of $\Proj^*A$. The inverse map to this correspondence is induced by
the composite
   $$\Proj^*A\supseteq U\bl\phi\mapsto\cc Q_{ U}
     \bl\zeta\mapsto\cc S(\cc Q_{ U}):=\{M\in\qgr A \mid(M)_q\subseteq\cc Q_{ U}\}$$
where $\zeta$ yields the inverse to the correspondence induced by
$\tau$ and $\psi$ yields the inverse to the correspondence induced
by $\phi$ (see Lemma~\ref{lll}).
\end{proof}

To conclude the section, we give a classification of the tensor
localizing subcategories of finite type in $\QGr A$ (=tensor Serre
subcategories $\cc L$ closed under direct limits such that $\QGr
A/\cc L\to\QGr A$ respects direct limits) in terms of open sets of
$\Proj^*A$. For commutative (ungraded) noetherian rings, $R$, a
similar classification of the localizing subcategories in $\Rfp$ in
terms of open sets of $\spec^*R$ has been shown by Hovey \cite{Ho}
and generalized to commutative (ungraded) coherent rings by Garkusha
and Prest~\cite[2.4]{GP}.

\begin{cor}\label{fff}%%\marginpar{fff}
The assignments
   $$\QGr A\supseteq\cc L\mapsto\bigcup_{M\in\cc L}\supp_A(M)$$
and
   $$\Proj^*A\supseteq U\mapsto
     \{\lp_\lambda M_\lambda\mid M_\lambda\in\qgr A,\supp_A(M_\lambda)\subseteq U\}$$
induce bijections between
\begin{itemize}
\item{$\diamond$} the set of all tensor localizing subcategories of
finite type in $\QGr A$,

\item{$\diamond$} the set of all open subsets $
U\subseteq\Proj^*A$.
\end{itemize}
\end{cor}

\begin{proof}
By~\cite[2.8]{H} and~\cite[2.10]{Kr1} there is a 1-1 correspondence
between the tensor Serre subcategories of $\qgr A$ and the
localizing subcategories of finite type in $\QGr A$. This
correspondence is given by
   $$\cc S\mapsto\vec{\cc S}:=\{\lp_\lambda M_\lambda\mid M_\lambda\in\cc
   S\}\ \ \ \textrm{ and }\ \ \ \cc L\mapsto\cc L\cap\qgr A.$$
Since the functor of $P$-localization with $P\in\Proj A$ commutes
with direct limits, we see that
   $$\bigcup_{M\in\cc L}\supp_A(M)=\bigcup_{M\in\cc L\cap\qgr A}\supp_A(M)$$ Now our assertion
follows from Theorem~\ref{tol}.
\end{proof}

\section{The prime spectrum of an ideal lattice}

Inspired by recent work of Balmer~\cite{B}, Buan, Krause, and
Solberg~\cite{BKS} introduce the notion of an ideal lattice and
study its prime ideal spectrum. Applications arise from studying
abelian or triangulated tensor categories.

\begin{defs}[Buan, Krause, Solberg~\cite{BKS}] {\rm
An {\it ideal lattice\/} is by definition a partially ordered set
$L=(L,\leq)$, together with an associative multiplication $L\times
L\to L$, such that the following holds.
\begin{enumerate}
\item[(L1)] The poset $L$ is a {\it complete lattice\/}, that is,
$$\sup A = \bigvee_{a\in A} a\quad\text{and}\quad \inf A = \bigwedge_{a\in A}
a$$ exist in $L$ for every subset $A\subseteq L$.
\item[(L2)] The lattice $L$ is {\it compactly generated\/}, that is,
every element in $L$ is the supremum of compact elements.  (An
element $a\in L$ is {\em compact}, if for all $A\subseteq L$ with
$a\leq \sup A$ there exists some finite $A'\subseteq A$ with
$a\leq\sup A'$.)
\item[(L3)] We have for all $a,b,c\in L$
$$a(b\vee c)=ab\vee ac\quad\text{and}\quad (a\vee b)c=ac\vee bc.$$
\item[(L4)] The element $1=\sup L$ is compact, and $1a=a=a1$ for all $a\in L$.
\item[(L5)] The product of two compact elements is again compact.
\end{enumerate}
A {\it morphism\/} $\phi\colon L\to L'$ of ideal lattices is a map
satisfying
\begin{gather*}\label{eq:mor}
\phi(\bigvee_{a\in A}a)=\bigvee_{a\in A}\phi(a)\quad \text{for}\quad
A\subseteq L, \\
\phi(1)=1\quad\text{and}\quad\phi(ab)=\phi(a)\phi(b)\quad\text{for}\quad
a,b\in L.\notag
\end{gather*}
}\end{defs}

Let $L$ be an ideal lattice. Following~\cite{BKS} we define the
spectrum of prime elements in $L$. An element $p\neq 1$ in $L$ is
 {\em prime} if $ab\leq p$ implies $a\leq p$ or $b\leq p$ for
all $a,b\in L$. We denote by $\spec L$ the set of prime elements in
$L$ and define for each $a\in L$
   $$V(a)=\{p\in\spec L\mid a\leq p\}\quad\text{and}\quad D(a)=\{p\in\spec L\mid a\not\leq p\}.$$
The subsets of $\spec L$ of the form $V(a)$ are closed under forming
arbitrary intersections and finite unions.  More precisely,
   $$V(\bigvee_{i\in\Omega} a_i)=\bigcap_{i\in\Omega} V(a_i)\quad\text{and}\quad V(ab)=V(a)\cup V(b).$$
Thus we obtain the {\it Zariski topology\/} on $\spec L$ by
declaring a subset of $\spec L$ to be {\it closed\/} if it is of the
form $V(a)$ for some $a\in L$. The set $\spec L$ endowed with this
topology is called the {\it prime spectrum\/} of $L$.  Note that the
sets of the form $D(a)$ with compact $a\in L$ form a basis of open
sets. The prime spectrum $\spec L$ of an ideal lattice $L$ is
spectral~\cite[2.5]{BKS}.

There is a close relation between spectral spaces and ideal
lattices. Given a topological space $X$, we denote by $L_{\open}(X)$
the lattice of open subsets of $X$ and consider the multiplication
map
   $$L_{\open}(X)\times L_{\open}(
 X)\to L_{\open}(X),\quad (U,V)\mapsto UV=U\cap V.$$
The lattice $L_{\open}(X)$ is complete.

\begin{prop}[Buan, Krause, Solberg~\cite{BKS}]\label{pr:openlattice}
Let $X$ be a spectral space. Then $L_{\open}(X)$ is an ideal
lattice. Moreover, the map
   $$X\to\spec L_{\open}(X),\quad x \mapsto X\setminus \overline{\{x\}},$$
is a homeomorphism.
\end{prop}

This construction is referred to as the soberification of $X$ and
the above result can be seen as part of the Stone Duality Theorem
(see, for instance, \cite{Jo}).

We deduce from the classification of tensor Serre subcategories of
$\qgr A$ (Theorem~\ref{tol}) the following.

\begin{prop}\label{prr}
Let $A$ be a coherent graded ring which is finitely generated as
an $A_0$-algebra. Then $L_{\serre}(\qgr A)$ is an ideal lattice.
\end{prop}

\begin{proof}
The space $\Proj A$ is spectral by Proposition~\ref{d3}. Thus
$\Proj^*A$ is spectral and $L_{\open}(\Proj^*A)$ is an ideal lattice
by Proposition~\ref{pr:openlattice}. By Theorem~\ref{tol} we have an
isomorphism $L_{\open}(\Proj^*A)\cong L_{\serre}(\qgr A)$, and
therefore $L_{\serre}(\qgr A)$ is an ideal lattice.
\end{proof}

\begin{cor}\label{prrco} The points of $\spec L_{\serre}(\qgr A)$ are the
$\cap$-irreducible tensor Serre subcategories of $\qgr A$ and the
map
   \begin{equation}\label{ggg}
    f:\Proj^*A\lra{}\spec L_{\serre}(\qgr A),\quad P\longmapsto\cc S_P=\{M\in\qgr A\mid M_P=0\}
   \end{equation}
is a homeomorphism of spaces.
\end{cor}

\begin{proof}
This is a consequence of Theorem~\ref{tol} and
Propositions~\ref{pr:openlattice} and~\ref{prr}.
\end{proof}

We write $\spec(\qgr A):=\spec^*L_{\serre}(\qgr A)$ and
$\supp(M):=\{\cc P\in\spec(\qgr A)\mid M\not\in\cc P\}$ for
$M\in\qgr A$. It follows from Corollary~\ref{prrco} that
   $$\supp_A(M)=f^{-1}(\supp(M)).$$
Following \cite{B,BKS}, we define a structure sheaf on $\spec(\qgr
A)$ as follows.  For an open subset $U\subseteq \spec(\qgr A)$,
let
   $$\cc S_U=\{M\in\qgr A\mid\supp(M)\cap U=\emptyset\}$$
and observe that $\cc S_U$ is a tensor Serre subcategory.  We
obtain a presheaf of rings on $\spec(\qgr A)$ by
   $$U\mapsto\End_{\qgr A/\cc S_U}(\cc O).$$ If $V\subseteq U$ are open subsets, then the restriction map
   $$\End_{\qgr A/\cc S_U}(\cc O)\to\End_{\qgr A/\cc S_V}(\cc O)$$
is induced by the quotient functor $\qgr A/\cc S_U\to\qgr A/\cc
S_V$. The sheafification is called the {\it structure sheaf\/} of
$\qgr A$ and is denoted by $\cc O_{\qgr A}$. This is a sheaf of
commutative rings by~\cite[1.7]{Su}. Next we have
by~\cite[7.1]{BKS}
   $$\cc O_{\qgr A,\cc P}\cong\End_{\qgr A/\cc P}(\cc O)\quad\text{for each}\quad\cc P\in\spec(\qgr A).$$

The next theorem says that the abelian category $\qgr A$ contains
all the necessary information to reconstruct the projective scheme
$(\Proj A,\cc O_{\Proj A})$.

\begin{thm}[Reconstruction]\label{coh}
Let $A$ be a coherent graded ring which is finitely generated as an
$A_0$-algebra. The map in~\eqref{ggg} induces an isomorphism
   $$f:(\Proj A,\cc O_{\Proj A})\lra{\sim}(\spec(\qgr A),\cc O_{\qgr A})$$
of ringed spaces.
\end{thm}

\begin{proof}
Fix an open subset $U\subseteq\spec(\qgr A)$ and consider the
composition of the functors
   $$F:\qgr A\to\Gr A\xrightarrow{\wt{(-)}}\Qcoh\Proj
     A\xrightarrow{(-)|_{f^{-1}(U)}}\Qcoh f^{-1}(U).$$
Here we denote, for any graded $A$-module $M$, by $\wt M$ its
associated sheaf. By definition, the stalk of $\wt M$ at a
homogeneous prime $P$ equals the degree zero part $M_{(P)}$ of the
localized module $M_P$. We claim that $F$ annihilates $\cc S_U$. In
fact, $M\in\cc S_U$ implies $f^{-1}(\supp(M))\cap
f^{-1}(U)=\emptyset$ and therefore $\supp_A(M)\cap
f^{-1}(U)=\emptyset$. Thus $M_{(P)}=0$ for all $P\in f^{-1}(U)$ and
therefore $F(M)=0$. It follows that $F$ factors through $\qgr A/\cc
S_U$ and induces a map $\End_{\qgr A/\cc S_U}(\cc O)\to\cc O_{\Proj
A}(f^{-1}(U))$ which extends to a map $\cc O_{\qgr A}(U)\to\cc
O_{\Proj A}(f^{-1}(U))$. This yields the morphism of sheaves
$f^\sharp\colon\cc O_{\qgr A}\to f_*\cc O_{\Proj A}$.

Now fix a point $P\in\Proj A$. Then $f^\sharp$ induces a map
$f^\sharp_P\colon\cc O_{\qgr A,f(P)}\to\cc O_{\Proj A,P}$. We have
an isomorphism $\cc O_{\qgr A,f(P)}\cong\End_{\qgr A/f(P)}(\cc O)$.
By construction,
   $$f(P)=\{M\in\qgr A\mid M_P=0\}.$$
One has an isomorphism
   $$\End_{\qgr A/f(P)}(\cc O)\cong\homs_{A}(\cc O,\cc O)_{(P)}\cong
     \homs_{A}(A,A)_{(P)}=\cc O_{\Proj A,P}.$$
We have used here Corollary~\ref{g111}. We conclude that
$f^\sharp_P$ is an isomorphism. It follows that $f$ is an
isomorphism of ringed spaces if the map $f:\Proj A\to\spec (\qgr
A)$ is a homeomorphism. This last condition is a consequence of
Theorem~\ref{tol} and Propositions~\ref{pr:openlattice}
and~\ref{prr}.
\end{proof}

To conclude the paper observe that one can similarly reconstruct
the affine scheme $(\spec R,\cc O_R)$ with $R$ a commutative
coherent ring out of Serre subcategories of the abelian category
$\rfp$ of finitely presented modules. Theorem~B in~\cite{GP}
implies that $L_{\serre}(\rfp)$ is an ideal lattice and that the
map
   $$f:\spec^*R\lra{}\spec L_{\serre}(\rfp),\quad P\longmapsto\cc
S_P=\{M\in\rfp\mid M_P=0\}$$ is a homeomorphism of spaces. Next, one
proves, similarly to Theorem~\ref{coh}, that the map $f$ induces an
isomorphism
   $$(\spec R,\cc O_R)\lra{\sim}(\spec(\rfp),\cc O_{\rfp})$$
of ringed spaces, where the definition of $(\spec(\rfp),\cc
O_{\rfp})$ is like that of $(\spec(\qgr A),\cc O_{\qgr A})$.

\end{document}